\documentclass[sn-mathphys]{sn-jnl}
\jyear{2023}

\theoremstyle{thmstyleone}

\theoremstyle{thmstyletwo}

\theoremstyle{thmstylethree}%

\usepackage{appendix}

\usepackage{graphicx}
\usepackage{hyperref}
\usepackage{subcaption}
\usepackage{bm}
\usepackage{xcolor}
\usepackage{bm,amsmath,amssymb,amsfonts,dsfont}
\usepackage[boxruled,ruled,vlined]{algorithm2e}
\DontPrintSemicolon

\newcommand{\variables}{x}
\newcommand{\Variables}{X}
\newcommand{\continuousvariables}{x^c}
\newcommand{\continuousvariablesset}{\mathcal{X}^c}
\newcommand{\Continuousvariables}{X^c}
\newcommand{\constraints}{g}

\newcommand{\catalogvariable}{c}
\newcommand{\catalog}{\mathcal{C}}
\newcommand{\properties}{x^p}
\newcommand{\propertiesset}{\mathcal{X}^p}
\newcommand{\Properties}{{X^p}}
\newcommand{\numbervariables}{n_c}
\newcommand{\numberproperties}{n_p}

\newcommand\equaldef{\stackrel{\text{def}}{=}}

\newcommand{\intervalbox}[1]{\bm{#1}}
\newcommand{\clutched}[1]{\hat{#1}}

\begin{document}

\title[Interval constraint programming for catalog-based categorical optimization]{Interval constraint programming for globally solving catalog-based categorical optimization}

\author*{\fnm{Charlie} \sur{Vanaret}}\email{vanaret@zib.de}
\affil*{\orgdiv{Department of Applied Algorithmic Intelligence Methods}, \orgname{Zuse-Institut Berlin}, \city{Berlin}, \postcode{14195}, \country{Germany} \\
ORCID: 0000-0002-1131-7631}

\artnote{Dedicated to Leila}

\abstract{
In this article, we propose an interval constraint programming method for globally solving catalog-based categorical optimization problems. It supports catalogs of arbitrary size and properties of arbitrary dimension, and does not require any modeling effort from the user.
A novel catalog-based contractor (or filtering operator) guarantees consistency between the categorical properties and the existing catalog items. This results in an intuitive and generic approach that is exact, rigorous (robust to roundoff errors) and can be easily implemented in an off-the-shelf interval-based continuous solver that interleaves branching and constraint propagation.  
We demonstrate the validity of the approach on a numerical problem in which a categorical variable is described by a two-dimensional property space.
A Julia prototype is available as open-source software under the MIT license at \url{https://github.com/cvanaret/CateGOrical.jl}~.
}


\maketitle

\section{Introduction}

\subsection{Motivation}

Catalog-based categorical optimization -- categorical variables are described by a tuple of properties that assume values within a discrete non-ordered catalog -- arises in numerous structural optimization applications: optimal truss design~\cite{gao2018categorical}, tyre selection~\cite{nedvelkova2019splitting}, aircraft design~\cite{barjhoux2017mixed}, product family design~\cite{lindroth2011pure}, optimization of thermal insulation systems~\cite{abhishek2010modeling}, etc. While problems with discrete variables are noticeably harder to solve than continuous problems, categorical variables add another level of complexity, since the core technique for solving discrete optimization problems -- relaxation of the integrality constraints -- does not hold.

\begin{figure}[htbp!]
\centering
\includegraphics[width=0.75\columnwidth]{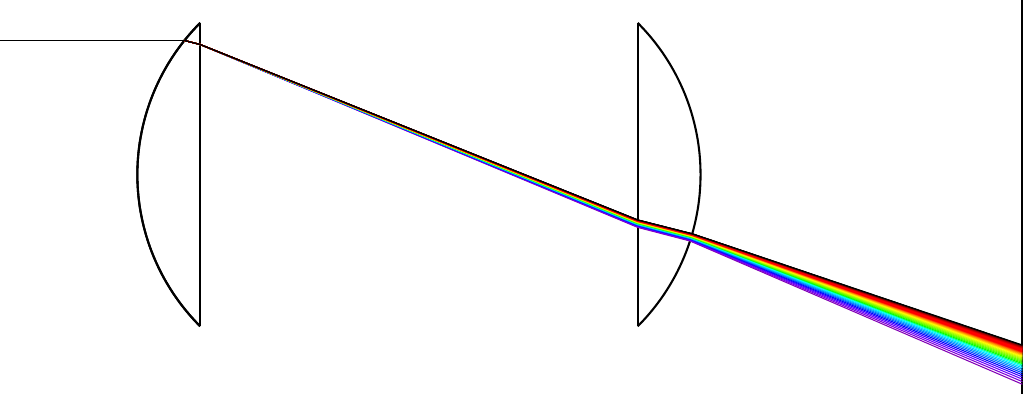}
\vspace{0.3cm}
\caption{Chromatic aberration in a refractive optical system.}
\label{fig:chromatic-aberration}
\end{figure}

A particular application in which categorical variables arise is the design of optical systems~\cite{agurok2019multi}: it consists in finding the parameters of a sequence of optical elements (e.g. their curvatures, distances and materials) such that the image of an object on the image plane minimizes the optical aberrations. Light propagation through the optical system is governed by the laws of reflection and refraction (Snell's laws): refracted light is bent at the interface between two media of different refractive indices (a number that describes how fast light travels through the material). The refractive index $n$ of a material is a function of the wavelength $\lambda$ of light -- a phenomenon known as \textit{dispersion} -- which causes different wavelengths (colors) to refract at different angles. A refractive optical system therefore exhibits \textit{chromatic aberration}, a failure to focus all colors at the same point. Figure~\ref{fig:chromatic-aberration} illustrates chromatic aberration for a system with two lenses: the incoming parallel rays with different wavelengths (the whole visible spectrum) are focused at different points on the image plane (the vertical line).

\begin{figure}[htbp!]
\centering
\includegraphics[width=0.6\columnwidth]{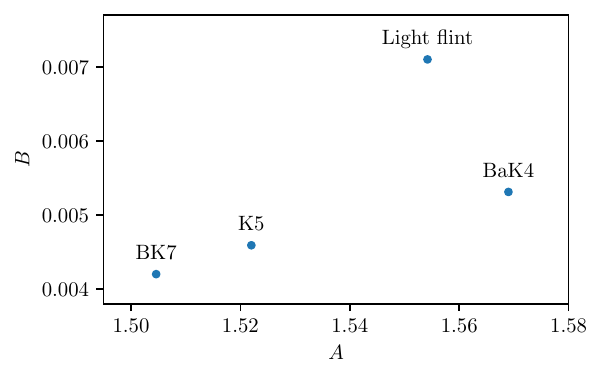}
\caption{Common optical glasses in the space of Cauchy coefficients.}
\label{fig:cauchy}
\end{figure}

Several (semi-)empirical models attempt to describe dispersion for a particular material. For instance, Cauchy's equation in its simplest form states:
\begin{equation*}
n(\lambda) = A + \frac{B}{\lambda^2},
\end{equation*}
where the coefficients $A$ and $B$ (in $\mu m^2$) are determined experimentally. Table~\ref{tab:cauchy-coefficients} lists the values of $A$ and $B$ for four common types of glasses~\cite{lin2017advanced}. Solving an optical design problem thus consists in picking optical materials (glasses) whose underlying properties (in our case, the Cauchy coefficients $A$ and $B$) appear in the equations that govern the propagation of light. These properties assume values within a finite set (the catalog) that, in the general case, cannot be ordered (Figure~\ref{fig:cauchy}).

\begin{table}[htbp!]
\caption{Cauchy coefficients of common optical glasses.}
\label{tab:cauchy-coefficients}
\centering
\begin{tabular}{lcc}
\hline
Material                   & $A$      & $B$ \\
\hline
Borosilicate glass (BK7)   & 1.5046   & 0.00420 \\
Hard crown glass (K5) 	   & 1.5220   & 0.00459 \\
Light flint                & 1.5542   & 0.00710 \\
Barium crown glass (BaK4)  & 1.5690   & 0.00531 \\
\hline
\end{tabular}
\end{table}

\subsection{Problem definition}

We consider the following mixed categorical optimization problem\footnote{For the sake of clarity and without loss of generality, we describe the problem with a single categorical variable and a single catalog. The proposed strategy however supports an arbitrary number of categorical variables and catalogs.}:
\begin{equation}
\label{eq:original-problem}
\begin{aligned}
& \underset{\continuousvariables, \catalogvariable}{\text{min}} & & f(\continuousvariables, \catalogvariable) & & \\
& \text{s.t.} 	& & \constraints(\continuousvariables, \catalogvariable) = 0 & & \\
&            	& & \catalogvariable \in \catalog & & \\
&               & & \continuousvariables \in \continuousvariablesset, & & \\
\end{aligned}
\end{equation}
where
\begin{itemize}
\item $\continuousvariables \in \continuousvariablesset \subset \mathbb{R}^{\numbervariables}$ is a vector of $\numbervariables$ continuous optimization variables;
\item $\catalogvariable$ is a categorical optimization variable;
\item $f$ is the (possibly nonlinear) objective function;
\item $\constraints(\continuousvariables, \catalogvariable) = 0$ are (possibly nonlinear) constraints in vector form; 
\item $\catalogvariable \in \catalog$ is a catalog constraint: $c$ should be assigned an item from the finite catalog $\catalog$ of size $\lvert \catalog \rvert$ (Table~\ref{tab:catalog}). We assume that all catalog items are represented by the same set of properties $\properties \in \propertiesset \subset \mathbb{R}^{\numberproperties}$.
\end{itemize}

\begin{table}[htbp!]
\caption{Catalog $\catalog$.}
\label{tab:catalog}
\centering
\begin{tabular}{c|ccc}
\hline
Item                       & Property $\properties_1$ & \ldots & Property $\properties_{\numberproperties}$ \\
\hline
1                          & value                    & \ldots & value \\
\vdots 	                  & \vdots                   & \vdots & \vdots \\
$\lvert \catalog \rvert$   & value                    & \ldots & value \\
\hline
\end{tabular}
\end{table}

Note that the catalog constraint $\catalogvariable \in \catalog$ is known as a \textit{table constraint} in the constraint programming community~\cite{lecoutre2012path}.

\subsection{State-of-the-art methods for categorical optimization}

\subsubsection{Heuristics}
Evolutionary and genetic algorithms have been among the first heuristics to automatize catalog-based optimization in industrial applications. While some authors adopted a generic approach~\cite{brown1993solving,carlson1996genetic}, others tailored the encoding (a hierarchical structure), crossover (a trial-and-error selection of breakpoints) and mutation operators to catalog design~\cite{carlson1998using}. 
\cite{lindroth2011pure} combined a local search and a global search for categorical optimization: a pattern search technique that iteratively improves the current point by visiting its categorical neighborhood is embedded within a discrete global descent framework.
In~\cite{gao2018categorical}, the authors interpreted the categorical properties as points in a high-dimensional metric space. They performed dimensionality reduction by using a manifold learning technique that preserves geodesic distances. They then defined new operators for genetic algorithms based on shortest paths in the resulting low-dimensional graph.
\cite{nedvelkova2019splitting} described a splitting strategy for simulation-based optimization that exploits convex relaxations of the discrete search space to evaluate approximate lower bounds. The bounds are then used in a branching heuristic within a multilevel coordinate search.

\subsubsection{Exact methods}
\label{sec:exact-methods}
Several families of methods solve categorical optimization to global optimality. \cite{abhishek2010modeling} modeled the mixed categorical problem with binary variables and replaced the assignments of categorical variables with convex combinations of the catalog entries. They solved the resulting mixed-integer nonlinear programming (MINLP) models with standard branch-and-bound methods via the AMPL modeling language.
\cite{lecoutre2012path} introduced STR3, a novel approach that built upon its predecessors. It enforces Generalized Arc Consistency on a positive table constraint (a catalog of permitted combinations of values) by maintaining a ``dual'' representation of the table constraint: each property value is assigned the list of items in which the value occurs. STR3 is optimal in the sense that it avoids unnecessary traversal of the table.
\cite{barjhoux2017mixed} introduced a branch-and-bound approach that instantiates a categorical variable at each level of the search tree. It evaluates lower bounds by relaxing the constraints that depend on undecided categorical variables and progressively inserts them back into the problem when the corresponding categorical variables are assigned a value.

Unfortunately, the immense majority of MINLP solvers compute nonrigorous relaxations based on floating-point arithmetic and are subject to roundoff errors~\cite{neumaier2005comparison}, which may result in erroneous (suboptimal, ``over''-optimal or infeasible) results. For instance, one of the most robust MINLP solvers BARON~\cite{sahinidis1996baron}  is not fully rigorous.
On the other hand, most rigorous global optimization solvers based on interval analysis~\cite{Moore1966Interval} -- such as Ibex~\cite{Trombettoni2011Inner}, IBBA~\cite{Ninin2010Reliable}, GlobSol~\cite{Kearfott1996Rigorous} and Charibde~\cite{Vanaret2015Hybridization} -- do not readily support discrete optimization. The COCONUT Environment~\cite{neumaier2004complete}, a collection of rigorous interval techniques for MINLP is still under development.

There are however a few notable exceptions: ILOG Solver~\cite{Puget1994ilogsolver}, a discrete solver for constraint satisfaction problems (CSP), was extended with the interval solver Numerica~\cite{van1997numerica} (not maintained).
Choco~\cite{jussien2008choco}, a mixed-integer constraint programming (CP) solver, was later interfaced with Ibex~\cite{chabert2009contractor} to benefit from the latest advances in continuous CSP techniques and reduce the costs of implementation and maintenance.
In both cases, a mature discrete solver served as a starting point and was later enhanced with continuous CP capabilities. As we describe in the next section, we adopt the opposite perspective: state-of-the-art continuous CP codes can be extended with a modest discrete technique (\textsc{Clutch}) with virtually no implementation costs, while maintaining rigorousness.

\subsection{Contributions and outline}

We propose a \emph{purely continuous} rigorous interval-based method for solving mixed categorical optimization problems that can be easily implemented in off-the-shelf rigorous solvers such as Ibex, IBBA, GlobSol and Charibde. The key idea of our approach is to notice that a categorical variable $c$ is described by a vector of properties, each of which lives in an ordered set. Instead of optimizing on $c$, we treat the properties $\properties$ themselves as continuous optimization variables and partition $\variables$ into continuous and property variables: $\variables = \begin{pmatrix} \continuousvariables \\ \properties \end{pmatrix}$. The corresponding optimization problem is given by:
\begin{equation}
\begin{aligned}
& \underset{\variables}{\text{min}} & & f(\continuousvariables, \properties) & & \\
& \text{s.t.} 	& & \constraints(\continuousvariables, \properties) = 0 & & \\
&            	& & \catalogvariable(\properties) \in \catalog & & \\
&               & & \variables = \begin{pmatrix} \continuousvariables \\ \properties \end{pmatrix} \text{ where } \continuousvariables \in \continuousvariablesset, \quad \properties \in \propertiesset. & & \\
\end{aligned}
\end{equation}

Our approach manipulates continuous relaxations of the properties $\properties$ and alternates between branching and constraint propagation phases on both $\continuousvariables$ and $\properties$. To ensure consistency between the properties $\properties$ and the existing catalog entries, we introduce a new filtering operator (or contractor~\cite{chabert2009contractor}), \textsc{Clutch} (Catalog LookUp Then Convex Hull), that performs a catalog lookup and reduces the current range of $\properties$ by considering the convex hull of the existing catalog items within this range.

We highlight the main advantages of our approach:
\begin{enumerate}
\item it is \textbf{exact}, \textbf{rigorous} and produces the \textbf{global optimum} within a given numerical tolerance $\varepsilon > 0$ even in the presence of roundoff errors;
\item \textbf{no modeling effort} (e.g. reformulation with binary variables~\cite{abhishek2010modeling}) is required from the user;
\item since we directly optimize the properties of the categorical variables, \textbf{no additional variable} is introduced in the model;
\item using continuous relaxations of the catalog properties \textbf{reduces the overestimation} when interval enclosures of the functions are computed; a comparison between continuous relaxations and formulations with binary variables is given in Appendix~\ref{sec:overestimation}.
\item \textsc{Clutch} is intuitive and \textbf{easy to implement} in off-the-shelf interval-based solvers that are otherwise not designed for discrete optimization.
\end{enumerate}

Section~\ref{sec:clutch} describes the novel catalog-based contractor \textsc{Clutch}. In Section~\ref{sec:icp-framework}, we show how \textsc{Clutch} naturally fits within an interval-based solver that interleaves branching and constraint propagation phases. In Section~\ref{sec:results}, we demonstrate the validity of our approach by solving two scenarios of a toy problem.

\section{\textsc{Clutch}: a catalog-based contractor}
\label{sec:clutch}

\textsc{Clutch} (Catalog LookUp Then Convex Hull) is a novel catalog-based contractor that enforces local consistency with respect to each individual catalog constraint $\catalogvariable(\properties) \in \catalog$.
It should be invoked whenever the range of $\properties$ is reduced during the arborescent exploration (upon branching and filtering) to ensure that the range of $\properties$ is consistent with existing catalog items. This ensures that subspaces of the search space that violate the catalog constraint $\catalogvariable(\properties) \in \catalog$ are discarded by refutation.

In the following, we call \textit{box} a Cartesian product (or vector) of intervals and use upper case symbols.
Formally, given a box $\intervalbox{\Properties}$ of categorical properties and a catalog $\catalog$, \textsc{Clutch} returns a box $\clutched{\intervalbox{\Properties}}$ that satisfies the \textit{correction property}:
\begin{equation}
\intervalbox{\Properties} \cap \catalog \quad \subset \quad \clutched{\intervalbox{\Properties}} \quad \subset \quad \intervalbox{\Properties}
\end{equation}
where $\intervalbox{\Properties} \cap \catalog$ are the catalog items that are in $\intervalbox{\Properties}$.

Since interval-based solvers manipulate compact boxes whose faces are parallel to the axes, \textsc{Clutch} returns a continuous relaxation of $\intervalbox{\Properties} ~\cap~ \catalog$ that is \emph{the smallest box that contains all items of $\intervalbox{\Properties} \cap \catalog$}:
\begin{equation}
\clutched{\intervalbox{\Properties}} := \square(\intervalbox{\Properties} \cap \catalog),
\end{equation}
where $\square$ is the convex hull operation. The algorithmic procedure is described in Algorithm~\ref{alg:clutch}.

\begin{algorithm}[htbp!]
\caption{\textsc{Clutch} contractor.}
\label{alg:clutch}
\SetAlgoVlined
\footnotesize
\KwIn{box $\intervalbox{\Properties}$, catalog $\catalog$}
\KwOut{filtered box $\clutched{\intervalbox{\Properties}}$}
Initialize $\clutched{\intervalbox{\Properties}} \gets \varnothing$ \;
\ForEach{catalog item with properties $\properties$ in $\catalog$} {
	\If{$\properties$ is in $\intervalbox{\Properties}$} {
        Compute the convex hull: $\clutched{\intervalbox{\Properties}} \gets \square(\clutched{\intervalbox{\Properties}}, \{\properties\})$
    }
} 
\Return $\clutched{\intervalbox{\Properties}}$ \;
\end{algorithm}

\textsc{Clutch} enforces hull (2B) consistency~\cite{Benhamou1999Revising} of the catalog constraint, that is the property of arc consistency for each bound of the variables that occur in a constraint. This is a consequence of the endpoint representation of intervals in interval arithmetic.
2B is a weaker consistency than Generalized Arc Consistency enforced by the STR methods (mentioned in Section~\ref{sec:exact-methods}). For this reason, STR has more filtering power than \textsc{Clutch}, albeit with higher computational cost and memory footprint.

Upon completion of \textsc{Clutch}, one of three situations will occur:
\begin{enumerate}
\item $\clutched{\intervalbox{\Properties}}$ is empty, that is there are no catalog items that belong to $\intervalbox{\Properties}$: $\intervalbox{\Properties}$ is inconsistent with the catalog constraint and is discarded;
\item $\clutched{\intervalbox{\Properties}}$ contains a single catalog item and is therefore exactly $\square(\intervalbox{\Properties} \cap \catalog)$ without overestimation: the categorical variable is assigned to this item;
\item $\clutched{\intervalbox{\Properties}}$ contains at least two distinct items: $\square(\intervalbox{\Properties} \cap \catalog)$ is therefore a relaxation of $\intervalbox{\Properties} \cap \catalog$ and further filtering or branching of the subspace is required.
\end{enumerate}

Note that to reduce the cost of catalog lookup, it is possible to maintain a list $L$ of non-increasing size that contains the catalog items within $\intervalbox{\Properties}$. Whenever $\intervalbox{\Properties}$ is contracted (upon branching or filtering), \textsc{Clutch} browses through $L$ instead of the whole catalog, which reduces the number of comparisons from $O(\lvert \catalog \rvert)$ (where $\lvert \catalog \rvert$ is the number of catalog items) to $O(\lvert L \rvert)$. The catalog items not contained in the contracted box are then discarded from $L$.

\section{Interval constraint programming for globally solving catalog-based categorical optimization}
\label{sec:icp-framework}

In this section, we show how \textsc{Clutch} naturally fits within an interval branch-and-contract framework described in a generic fashion in Algorithm~\ref{alg:framework}.
The generic operations are described in the following subsections: exploration of the search space, upper bounding, branching, lower bounding, and filtering.

\begin{algorithm}[htbp!]
\caption{Framework for mixed categorical optimization.}
\label{alg:framework}
\SetAlgoVlined
\footnotesize
\KwIn{domain $\intervalbox{D}$, objective $f$, constraints $\{\constraints(\variables) = 0, \catalogvariable(\properties) \in \catalog\}$, tolerance $\varepsilon > 0$}
\KwOut{feasible point $\tilde{x}$ such that $\overline{f(\tilde{x})} - f^* \le \varepsilon$}
\textsc{Filtering}: contract initial domain $\intervalbox{D}$ wrt constraints \;
\smallskip

\textsc{Lower bounding}: compute lower bound $lb$ of problem on $\intervalbox{D}$ \;
\smallskip

Set best known solution and upper bound $(\tilde{x}, \tilde{f}) \gets (\varnothing, +\infty)$ \;
Initialize queue $\mathcal{L} \gets \{(\intervalbox{D}, lb)\}$ \;
\Repeat{$\mathcal{L}$ empty} {
	 \textsc{Exploration}: extract $(\intervalbox{\Variables}, lb)$ from $\mathcal{L}$ \;
    \smallskip
    
    \textsc{Pruning}: \lIf{$\tilde{f} - \varepsilon < lb$} {
         discard $\intervalbox{X}$
      }
    \smallskip
    

    \textsc{Upper bounding}: find feasible point $x \in \intervalbox{X}$ \;
    \lIf{$\overline{f(x)} < \tilde{f} - \varepsilon$} {
         update best known solution: $(\tilde{x}, \tilde{f}) \gets (x, \overline{f(x)})$
      }
    \smallskip
    
    \textsc{Branching}: bisect $\intervalbox{X}$ into $\intervalbox{X}_1$ and $\intervalbox{X}_2$ \;
    \For{$i \in \{1, 2\}$} {
        \textsc{Filtering}: contract $\intervalbox{X}_i$ wrt $\constraints(\variables) = 0$ (HC4Revise, Fig.~\ref{fig:hc4}) and $\catalogvariable(\properties) \in \catalog$ (\textsc{Clutch}, Alg.~\ref{alg:clutch}) within a fixed-point algorithm (Alg.~\ref{alg:fixed-point}) \;
         \lIf{$\intervalbox{X}_i$ is empty} {
            discard $\intervalbox{X}_i$
         }
         \smallskip
         
         \textsc{Lower bounding}: compute lower bound $lb$ of problem on $\intervalbox{X}_i$ \;
         \smallskip
         
         \textsc{Pruning}: \lIf{$\tilde{f} - \varepsilon < lb$} {
            discard $\intervalbox{X}_i$
         }
         \smallskip
        
        Insert $(\intervalbox{X}_i, lb)$ into $\mathcal{L}$ \;
    }
} 
\Return $(\tilde{x}, \tilde{f})$ \;
\end{algorithm}

The global optimization problem is usually transformed into a numerical CSP that keeps track of the best known upper bound $\tilde{f}$ of the global minimum $f^*$.
A solution is sought within a tolerance $\varepsilon > 0$ of the global minimum, which can be enforced on each subspace with the dynamic constraint $f(\variables) \le \tilde{f} - \varepsilon$.
Consequently, points whose objective values do not improve upon $\tilde{f}$ by at least $\varepsilon$ are discarded. A corollary is that the box that contains the global minimizer might be discarded during the process. However, the algorithm is guaranteed to return an upper bound of $\tilde{f}$ within $\varepsilon$, together with the associated feasible point.

A Julia prototype of the described interval-based method for solving mixed categorical problems to global optimality is available as open-source software under the MIT license at \url{https://github.com/cvanaret/CateGOrical.jl}~.

\subsection{Exploration strategy}

The priority queue $\mathcal{L}$ is usually implemented as a binary heap, that is a complete binary tree in which the priority of a node is greater than the priority of its children. Inserting an element or extracting the element with the highest priority is carried out in logarithmic time. The choice of priority determines the exploration strategy of the search space (depth-first search, best-first search, breadth-first search, MaxDist~\cite{Vanaret2015Hybridization}, etc.)

\subsection{Upper bounding}

The objective value of any feasible point of the problem (provided that the feasible set is not empty) is an upper bound of the global minimum. The best known upper bound $\tilde{f}$ of the global minimum can be improved by local search algorithms or metaheuristics in an integrated~\cite{zhang2007new} or cooperative manner~\cite{gallardo2007hybridization,blum2011hybrid,cotta1995hybridizing,alliot2012finding,Vanaret2015Hybridization}.
A cheaper strategy consists in systematically updating $\tilde{f}$ with the objective value of the midpoint of the current subspace, provided that it is feasible~\cite{ichida1979interval}.

We implemented a simple and cheap strategy to generate a possible feasible point $x = \begin{pmatrix} \continuousvariables \\ \properties \end{pmatrix}$ within a given box $\intervalbox{\Variables} = \begin{pmatrix} \intervalbox{\Continuousvariables} \\ \intervalbox{\Properties} \end{pmatrix}$. The components $\properties$ are set to some catalog item within $\intervalbox{\Properties}$ (there should exist one, otherwise $\intervalbox{\Variables}$ would have been pruned by \textsc{Clutch}). A round of constraint propagation with respect to the general constraints tightens the domains of the continuous components $\intervalbox{\Continuousvariables}$. This smaller box is by construction always catalog-feasible. Its midpoint, provided that it is feasible, can then be used to update $\tilde{f}$.

\subsection{Branching}

Branching consists in partitioning the current subspace into several (usually two) subspaces. There exist a variety of strategies, among which:
\begin{itemize}
\item the variable with the largest domain is partitioned;
\item the variables are partitioned one after the other in a round-robin fashion;
\item the smear heuristic picks the variable that has the largest estimated impact (a measure based on the interval gradient) on the objective function~\cite{kearfott1990algorithm};
\item the variable $\variables$ that leads to the highest relaxation error of a nonlinear term featuring $\variables$ is partitioned.
\end{itemize}
In the proposed strategy, branching is carried out on both continuous variables $\continuousvariables$ and properties $\properties$.

\subsection{Lower bounding}

The most commonly used technique to evaluate a lower bound of a constrained problem on a given subspace is the Reformulation Linearization Technique (RLT)~\cite{sherali2013reformulation}. It consists in generating a convex outer approximation (by linearization or convexification) of the problem, then solving the convex relaxed problem. Its interval counterpart X-Newton~\cite{Araya2012Contractor} implements rigorous linearizations using interval arithmetic~\cite{sunaga1958theory} and a cheap duality-based postprocessing step~\cite{Neumaier2004Safe} for evaluating rigorous lower bounds even in the presence of roundoff errors.

A much cheaper approach, albeit of coarser quality, is to compute an enclosure of the range of the sole objective function on a given subspace using interval arithmetic. It amounts to computing a lower bound of the unconstrained problem and may therefore largely overestimate the actual range.

\subsection{Filtering}
\label{sec:filtering}

Filtering algorithms (or contractors) reduce the domains of the variables with respect to individual constraints (local consistency) or all the constraints simultaneously (global consistency) by discarding inconsistent values. They usually stem from the mixed integer, constraint programming and numerical analysis communities, and may exploit the syntax tree, range or monotonicity of an individual constraint, or perform a convexification of the whole system.

\begin{figure}[htbp!]
\centering
\begin{subfigure}[b]{\textwidth}
\centering
\footnotesize
\def\svgwidth{0.65\columnwidth}
\begingroup%
  \makeatletter%
  \providecommand\color[2][]{%
    \errmessage{(Inkscape) Color is used for the text in Inkscape, but the package 'color.sty' is not loaded}%
    \renewcommand\color[2][]{}%
  }%
  \providecommand\transparent[1]{%
    \errmessage{(Inkscape) Transparency is used (non-zero) for the text in Inkscape, but the package 'transparent.sty' is not loaded}%
    \renewcommand\transparent[1]{}%
  }%
  \providecommand\rotatebox[2]{#2}%
  \ifx\svgwidth\undefined%
    \setlength{\unitlength}{509.37138672bp}%
    \ifx\svgscale\undefined%
      \relax%
    \else%
      \setlength{\unitlength}{\unitlength * \real{\svgscale}}%
    \fi%
  \else%
    \setlength{\unitlength}{\svgwidth}%
  \fi%
  \global\let\svgwidth\undefined%
  \global\let\svgscale\undefined%
  \makeatother%
  \begin{picture}(1,0.38108912)%
    \put(0,0){\includegraphics[width=\unitlength]{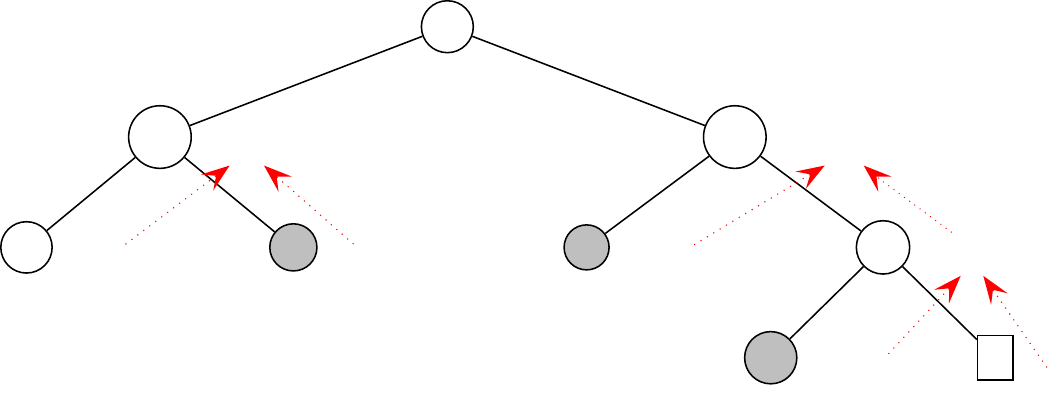}}%
    \put(0.01732619,0.13500928){\makebox(0,0)[lb]{\smash{2}}}%
    \put(0.92804382,0.03048191){\makebox(0,0)[lb]{\smash{2}}}%
    \put(0.08082718,0.12484952){\color[rgb]{0,0,0}\makebox(0,0)[lb]{\smash{$[2, 2]$}}}%
    \put(0.30542058,0.12205594){\color[rgb]{0,0,0}\makebox(0,0)[lb]{\smash{$[0, 20]$}}}%
    \put(0.19376295,0.2387516){\color[rgb]{0,0,0}\makebox(0,0)[lb]{\smash{$[0, 40]$}}}%
    \put(0.77397913,0.01423112){\color[rgb]{0,0,0}\makebox(0,0)[lb]{\smash{$[-10, 10]$}}}%
    \put(0.86734463,0.12968439){\color[rgb]{0,0,0}\makebox(0,0)[lb]{\smash{$[0, 100]$}}}%
    \put(0.61025045,0.12431424){\color[rgb]{0,0,0}\makebox(0,0)[lb]{\smash{$[0, 16]$}}}%
    \put(0.74120402,0.23576409){\color[rgb]{0,0,0}\makebox(0,0)[lb]{\smash{$[-100, 16]$}}}%
    \put(0.54352708,0.13910705){\color[rgb]{0,0,0}\makebox(0,0)[lb]{\smash{$z$}}}%
    \put(0.26643484,0.13966798){\color[rgb]{0,0,0}\makebox(0,0)[lb]{\smash{$x$}}}%
    \put(0.71797178,0.03645952){\color[rgb]{0,0,0}\makebox(0,0)[lb]{\smash{$y$}}}%
    \put(0.13764866,0.24141803){\color[rgb]{0,0,0}\makebox(0,0)[lb]{\smash{$\times$}}}%
    \put(0.8202828,0.14179946){\color[rgb]{0,0,0}\makebox(0,0)[lb]{\smash{$\wedge$}}}%
    \put(0.67870775,0.24175459){\color[rgb]{0,0,0}\makebox(0,0)[lb]{\smash{$-$}}}%
    \put(0.40834645,0.34664577){\color[rgb]{0,0,0}\makebox(0,0)[lb]{\smash{$=$}}}%
    \put(0.9850797,0){\color[rgb]{0,0,0}\makebox(0,0)[lb]{\smash{2}}}%
  \end{picture}%
\endgroup%
\caption{Bottom-up evaluation phase.}
\end{subfigure}
\par\bigskip
\begin{subfigure}[b]{\textwidth}
\centering
\footnotesize
\def\svgwidth{0.65\columnwidth}
\begingroup%
  \makeatletter%
  \providecommand\color[2][]{%
    \errmessage{(Inkscape) Color is used for the text in Inkscape, but the package 'color.sty' is not loaded}%
    \renewcommand\color[2][]{}%
  }%
  \providecommand\transparent[1]{%
    \errmessage{(Inkscape) Transparency is used (non-zero) for the text in Inkscape, but the package 'transparent.sty' is not loaded}%
    \renewcommand\transparent[1]{}%
  }%
  \providecommand\rotatebox[2]{#2}%
  \ifx\svgwidth\undefined%
    \setlength{\unitlength}{362.13999023bp}%
    \ifx\svgscale\undefined%
      \relax%
    \else%
      \setlength{\unitlength}{\unitlength * \real{\svgscale}}%
    \fi%
  \else%
    \setlength{\unitlength}{\svgwidth}%
  \fi%
  \global\let\svgwidth\undefined%
  \global\let\svgscale\undefined%
  \makeatother%
  \begin{picture}(1,0.40857711)%
    \put(0,0){\includegraphics[width=\unitlength]{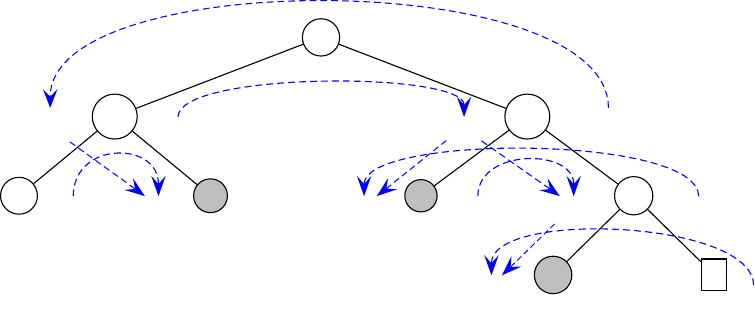}}%
    \put(0.01742878,0.13657604){\makebox(0,0)[lb]{\smash{2}}}%
    \put(0.2686391,0.13924958){\makebox(0,0)[lb]{\smash{$x$}}}%
    \put(0.54874327,0.13924958){\makebox(0,0)[lb]{\smash{$z$}}}%
    \put(0.72425298,0.03650867){\makebox(0,0)[lb]{\smash{$y$}}}%
    \put(0.93844251,0.02938636){\makebox(0,0)[lb]{\smash{2}}}%
    \put(0.06450589,0.12094325){\color[rgb]{0,0,0}\makebox(0,0)[lb]{\smash{$[2, 2]$}}}%
    \put(0.82825407,0.1418322){\color[rgb]{0,0,0}\makebox(0,0)[lb]{\smash{$\wedge$}}}%
    \put(0.68495814,0.24259323){\color[rgb]{0,0,0}\makebox(0,0)[lb]{\smash{$-$}}}%
    \put(0.41147991,0.34935065){\color[rgb]{0,0,0}\makebox(0,0)[lb]{\smash{$=$}}}%
    \put(0.13744118,0.24287348){\color[rgb]{0,0,0}\makebox(0,0)[lb]{\smash{$\times$}}}%
    \put(0.16028395,0.1213851){\color[rgb]{0,0,0}\makebox(0,0)[lb]{\smash{$\mathbf{[0, 8]}$}}}%
    \put(0.01674078,0.23545791){\color[rgb]{0,0,0}\makebox(0,0)[lb]{\smash{$\mathbf{[0, 16]}$}}}%
    \put(0.21293561,0.22111615){\color[rgb]{0,0,0}\makebox(0,0)[lb]{\smash{$[0, 40]$}}}%
    \put(0.99160548,0){\color[rgb]{0,0,0}\makebox(0,0)[lb]{\smash{2}}}%
    \put(0.59736913,0.01068302){\color[rgb]{0,0,0}\makebox(0,0)[lb]{\smash{$\mathbf{[-4, 4]}$}}}%
    \put(0.7639311,0.01121804){\color[rgb]{0,0,0}\makebox(0,0)[lb]{\smash{$[-10, 10]$}}}%
    \put(0.69995377,0.1213851){\color[rgb]{0,0,0}\makebox(0,0)[lb]{\smash{$\mathbf{[0, 16]}$}}}%
    \put(0.77038578,0.23331792){\color[rgb]{0,0,0}\makebox(0,0)[lb]{\smash{$[-100, 16]$}}}%
    \put(0.3125863,0.12094325){\color[rgb]{0,0,0}\makebox(0,0)[lb]{\smash{$[0, 20]$}}}%
    \put(0.43091071,0.1213851){\color[rgb]{0,0,0}\makebox(0,0)[lb]{\smash{$\mathbf{[0, 16]}$}}}%
    \put(0.59262303,0.12094325){\color[rgb]{0,0,0}\makebox(0,0)[lb]{\smash{$[0, 16]$}}}%
    \put(0.89076396,0.12094325){\color[rgb]{0,0,0}\makebox(0,0)[lb]{\smash{$[0, 100]$}}}%
    \put(0.55950394,0.2293312){\color[rgb]{0,0,0}\makebox(0,0)[lb]{\smash{$\mathbf{[0, 16]}$}}}%
  \end{picture}%
\endgroup%
\caption{Top-down propagation phase.}
\end{subfigure}
\caption{FBBT/HC4Revise on the constraint $2x = z - y^2$.}
\label{fig:hc4}
\end{figure}

One of the most intuitive approaches, called evaluation-propagation algorithm~\cite{messine1997methodes,messine2004deterministic}, FBBT (Feasibility-Based Bounds Tightening)~\cite{belotti2010feasibility} or HC4Revise~\cite{Benhamou1999Revising}, draws its inspiration from relational interval arithmetic~\cite{cleary1987logical}. It removes values that cannot satisfy a constraint by performing two traversals of its syntax tree, in which leaves represent variables and constants, and nodes elementary operations. 
Figure~\ref{fig:hc4} illustrates the approach on the constraint $2x = z - y^2$ with $x \in [0, 20]$, $y \in [-10, 10]$ and $z \in [0, 16]$.


In order to contract a box with respect to a system of constraints, a fixed-point algorithm (an idempotent propagation loop) named AC-3~\cite{mackworth1977consistency} handles the constraints \textit{individually} and propagates range reductions through the system (Algorithm~\ref{alg:fixed-point}). 

\begin{algorithm}[htbp!]
\small
\caption{Fixed-point algorithm.}
\label{alg:fixed-point}
\footnotesize
\KwIn{box $\intervalbox{\Variables}$, contractors $\mathcal{L} = \{C_1, \ldots, C_m\}$}
\KwOut{contracted box}
Initialize the set of awake contractors $\mathcal{A} \gets \{C_1, \ldots, C_m\}$ \;
\Repeat{$\mathcal{L}$ is empty} {
	Extract a contractor $C$ from $\mathcal{A}$ \;
	Filter the box: $\intervalbox{\Variables}' \gets C(\intervalbox{\Variables})$ \;
	\If{$\intervalbox{\Variables}'$ sufficiently smaller than $\intervalbox{\Variables}$} {
		Add to $\mathcal{A}$ the contractors of $\mathcal{L}$ that involve contracted variables of $\intervalbox{\Variables}'$ \;
		Update the box: $\intervalbox{\Variables} \gets \intervalbox{\Variables}'$ \;
	}
}
\Return $\intervalbox{\Variables}$ \;
\end{algorithm}

The extension of the proposed method to multiple categorical variables and catalogs is straightforward: the \textsc{Clutch} contractor is attached as a filtering procedure to each individual catalog constraint. The general constraints and the catalog constraints are then handled by the fixed-point algorithm.

Convexification techniques, called OBBT (Optimization-Based Bound Tightening)~\cite{zamora1999branch} or X-Newton~\cite{Araya2012Contractor}, handle the constraints \textit{simultaneously} by generating convex under- and overestimators of the constraints on a given subspace. The lower (resp. upper) bound of the variable $\variables_i$ is reduced by minimizing (resp. maximizing) $\variables_i$ over the convex feasible set. For a comprehensive survey, see~\cite{Vanaret2015Hybridization}.

\section{Numerical results}
\label{sec:results}

In this section, we describe a toy problem and solve it to global optimality using the proposed approach. Although the problem is trivial and can be solved by simple enumeration of the catalog items, the example demonstrates that continuous relaxations of the catalog properties naturally fit within an interval-based framework. We also include diagrams that provide a visual interpretation of the \textsc{Clutch} contractor.
We provide a comparison between our approach and the approach with binary decision variables described in~\cite{abhishek2010modeling}; the resulting mixed model is solved with SCIP~\cite{SCIP}, one of the most efficient open-source MINLP solvers.
For a fair comparison between both approaches, the SCIP presolve was disabled.

\subsection{Toy problem}
\label{sec:toy-problem}

We introduce the following toy problem:
\begin{equation}
\begin{aligned}
& \underset{x}{\text{min}} & & f(\variables) \equaldef (\properties_1)^3 & & \\
& \text{s.t.} 	& & \constraints(\variables) \equaldef \continuousvariables_1 - (\properties_2)^2 - 2 \properties_1 = 0 & & \\
&               & & \catalogvariable(\properties_1, \properties_2) \in \catalog & & \\
&               & & \continuousvariables_1 \in [0, 16], \properties_1 \in [0, 20], \properties_2 \in [-10, 10], 
\end{aligned}
\end{equation}
where $(\properties_1, \properties_2)$ are categorical properties and $\continuousvariables_1$ is a continuous variable. Note that the constraint $g$ is the same as in Figure~\ref{fig:hc4}.

We propose two scenarios, each with a different catalog:
\begin{itemize}
\item in \textbf{scenario~1}, the catalog $\catalog$ contains five items (Table~\ref{tab:scenario1-catalog}). The feasible set of the problem contains a single catalog item (item 2);
\item in \textbf{scenario~2}, the catalog $\catalog$ of the first scenario is augmented with an item described by properties $(\properties_1, \properties_2) = (1, -1)$ (Table~\ref{tab:scenario2-catalog}). The feasible set of the problem contains two catalog items (items 2 and 6).
\end{itemize}


\begin{minipage}{.5\linewidth}
   \tablecaptionfont
   \centering
    \captionof{table}{Scenario 1: catalog $\catalog$.}
    \label{tab:scenario1-catalog}
    \begin{tabular}{c|cc}
    \hline
    Item    & Property $\properties_1$  & Property $\properties_2$ \\
    \hline
    1       & 4             & -8 \\
    2       & 3             & 2 \\
    3       & 7             & -3 \\
    4       & 14            & 8 \\
    5       & 19            & -8 \\
    \hline
    \end{tabular}
    \vspace{0.4cm}
\end{minipage}%
\begin{minipage}{.5\linewidth}
   \tablecaptionfont
   \centering
    \captionof{table}{Scenario 2: catalog $\catalog$.}
    \label{tab:scenario2-catalog}
    \begin{tabular}{c|cc}
    \hline
    Item    & Property $\properties_1$  & Property $\properties_2$ \\
    \hline
    1       & 4             & -8 \\
    2       & 3             & 2 \\
    3       & 7             & -3 \\
    4       & 14            & 8 \\
    5       & 19            & -8 \\
    6       & 1             & -1 \\
    \hline
    \end{tabular}
\end{minipage}
\vspace{0.5cm}

The optimal solutions are:
\begin{itemize}
\item for scenario 1, $x^* = (10, 3, 2)$ and $f^* = f(x^*) = 27$;
\item for scenario 2, $x^* = (3, 1, -1)$ and $f^* = f(x^*) = 1$.
\end{itemize}

\begin{figure}[htbp!]
\centering
\begin{subfigure}[b]{0.49\textwidth}
\includegraphics[width=\columnwidth]{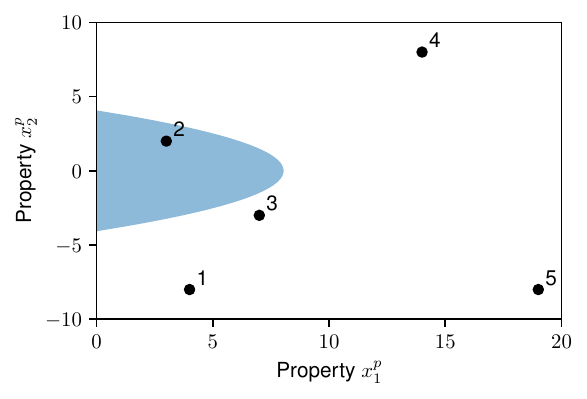}
\caption{Scenario 1.}
\label{fig:scenario1-properties}
\end{subfigure}
\begin{subfigure}[b]{0.49\textwidth}
\includegraphics[width=\columnwidth]{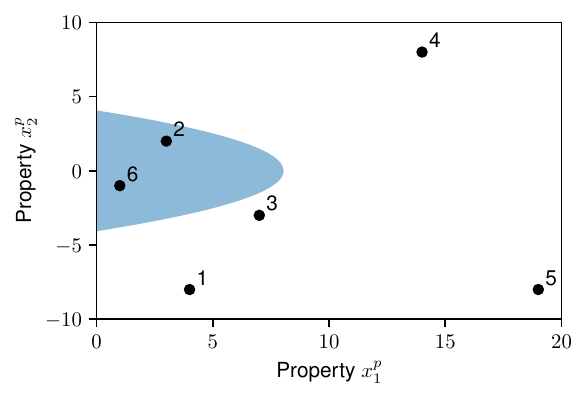}
\caption{Scenario 2.}
\label{fig:scenario2-properties}
\end{subfigure}
\caption{Property space of the toy problem.}
\label{fig:example-properties}
\end{figure}

The catalog items of scenarios 1 and 2 are represented in Figures~\ref{fig:scenario1-properties} and \ref{fig:scenario2-properties}, respectively. The blue domain is the feasible set of $\constraints$ (that is the set $\{ \continuousvariables_1 - (\properties_2)^2 - 2 \properties_1 = 0 \,\vert\, x \in [0, 16] \times [0, 20] \times [-10, 10] \}$) projected onto the $(\properties_1, \properties_2)$ plane.

In the following, we describe the sequence of steps carried out by the proposed approach to solve the two scenarios of the toy problem. Note that the steps depend on the heuristics used for each of the generic operations described in Algorithm~\ref{alg:framework}. Therefore, the Julia prototype \texttt{CateGOrical.jl} may not reproduce the exact sequence of steps described in this section.

\subsection{Scenario 1}
\label{sec:results-scenario1}

\subsubsection{Proposed approach}

\begin{figure}[htbp!]
\centering
\begin{subfigure}[b]{0.49\textwidth}
\includegraphics[width=\columnwidth]{example_scenario1}
\caption{Initial range of properties $\properties_1$ and $\properties_2$.}
\end{subfigure}
\hfill
\begin{subfigure}[b]{0.49\textwidth}
\includegraphics[width=\columnwidth]{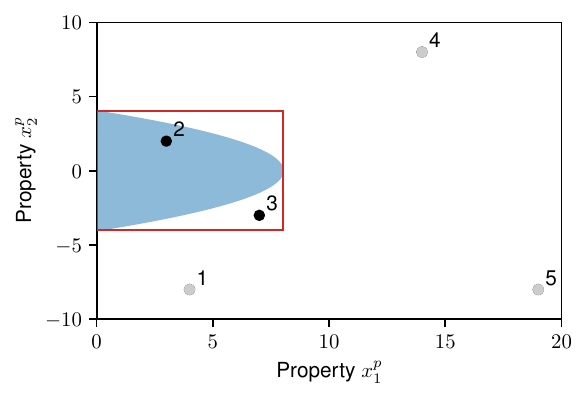}
\caption{HC4Revise with respect to $\constraints(\variables) = 0$.}
\end{subfigure} \\
\begin{subfigure}[b]{0.49\textwidth}
\includegraphics[width=\columnwidth]{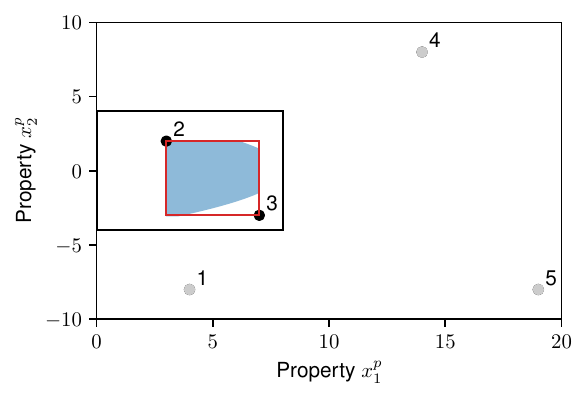}
\caption{\textsc{Clutch} with respect to $\catalogvariable(\properties) \in \catalog$.}
\end{subfigure}
\hfill
\begin{subfigure}[b]{0.49\textwidth}
\includegraphics[width=\columnwidth]{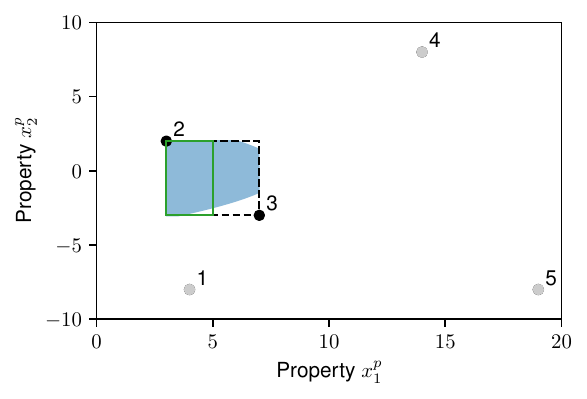}
\caption{Branching on $\properties_1$.}
\end{subfigure} \\
\begin{subfigure}[b]{0.49\textwidth}
\includegraphics[width=\columnwidth]{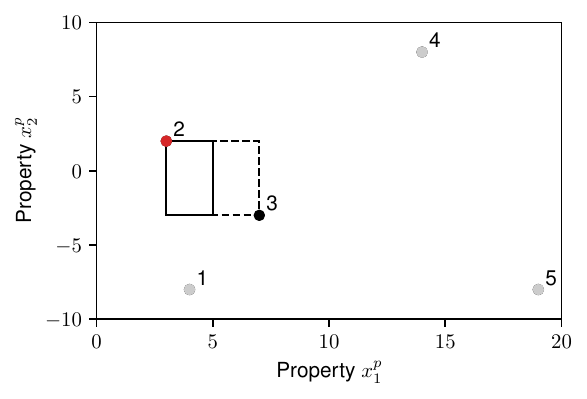}
\caption{\textsc{Clutch} with respect to $\catalogvariable(\properties) \in \catalog$.}
\end{subfigure}
\hfill
\begin{subfigure}[b]{0.49\textwidth}
\includegraphics[width=\columnwidth]{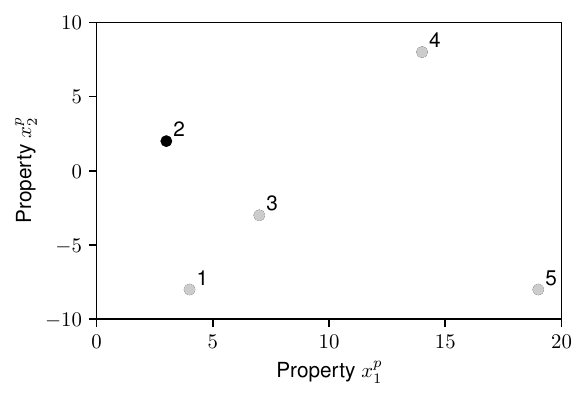}
\caption{HC4Revise with respect to $\constraints(\variables) = 0$.}
\end{subfigure}
\caption{Scenario 1: alternation of filtering and branching phases.}
\label{fig:scenario1}
\end{figure}

Table~\ref{tab:scenario1-phases} describes the successive filtering and branching phases (second column) that reduce the range of $x$ starting from $[0, 16] \times [0, 20] \times [-10, 10]$: the third column contains the current box for each phase and the fourth column presents the resulting contracted ranges. The first column indicates the corresponding subfigures in Figure~\ref{fig:scenario1}. The two subspaces resulting from the branching phase are highlighted in red.

\begin{table}[htbp!]
\caption{Scenario 1: alternation of filtering and branching phases.}
\label{tab:scenario1-phases}
\centering
\footnotesize
\begin{tabular}{clcc}
\hline
 & Phase                          & Current box                                & Result \\
\hline
(b) & HC4Revise($\constraints(\variables) = 0$)         & $[0, 16] \times [0, 20] \times [-10, 10]$  & $[0, 16] \times [0, 8] \times [-4, 4]$ \\
(c) & \textsc{Clutch}($\catalogvariable(\properties) \in \catalog$)   & $[0, 16] \times [0, 8] \times [-4, 4]$     & $[0, 16] \times [3, 7] \times [-3, 2]$ \\
    & HC4Revise($\constraints(\variables) = 0$)         & $[0, 16] \times [3, 7] \times [-3, 2]$     & $[6, 16] \times [3, 7] \times [-3, 2]$ \\
(d) & Branching on $\properties_1$       & $[6, 16] \times [3, 7] \times [-3, 2]$     & $[6, 16] \times \textcolor{red}{[3, 5]} \times [-3, 2] ~\cup$ \\
    &                          &                                            & $[6, 16] \times \textcolor{red}{[5, 7]} \times [-3, 2]$ \\
\hline
& \textsc{Clutch}($\catalogvariable(\properties) \in \catalog$)       & $[6, 16] \times \textcolor{red}{[5, 7]} \times [-3, 2]$     & $[6, 16] \times [7, 7] \times [-3, -3]$ \\
& HC4Revise($\constraints(\variables) = 0$)             & $[6, 16] \times [7, 7] \times [-3, -3]$    & $\varnothing$ \\
\hline
(e) & \textsc{Clutch}($\catalogvariable(\properties) \in \catalog$)   & $[6, 16] \times \textcolor{red}{[3, 5]} \times [-3, 2]$     & $[6, 16] \times [3, 3] \times [2, 2]$ \\
(f) & HC4Revise($\constraints(\variables) = 0$)         & $[6, 16] \times [3, 3] \times [2, 2]$      & $\mathbf{[10, 10] \times [3, 3] \times [2, 2]}$ \\
\hline
\end{tabular}
\end{table}

Figure~\ref{fig:scenario1} provides a visual interpretation of the filtering and branching phases. The current box is represented as a black rectangle, boxes contracted by HC4Revise or \textsc{Clutch} as red rectangles, boxes bisected during branching as green rectangles and discarded catalog items as gray points. As in Figure~\ref{fig:example-properties}, the feasible set is represented as a blue domain.

These results show that the first contraction step discards three of the five catalog items; convergence to the optimal (and only feasible) point ${x^* = (10, 3, 2)}$ is achieved by an alternation of HC4Revise (for $\constraints(x) = 0$) and \textsc{Clutch} (for $\catalogvariable(\properties) \in \catalog$) and a single branching step on $\properties_1$. The first subspace $[6, 16] \times [5, 7] \times [-3, 2]$ generated by branching is quickly discarded, while the second $[6, 16] \times [3, 5] \times [-3, 2]$ is subsequently contracted to the optimal solution. In this case, the assignment of the catalog variable to a catalog item does not require the use of the objective function.

\subsubsection{Model with binary decision variables}

The corresponding model with binary decision variables is given by:
\begin{equation}
\begin{aligned}
& \underset{\variables, b}{\text{min}} & & (\properties_1)^3 & & \\
& \text{s.t.} 	& & \continuousvariables_1 - (\properties_2)^2 - 2 \properties_1 = 0 & & \\
&               & & \properties_1 = 4 b_1 + 3 b_2 + 7 b_3 + 14 b_4 + 19 b_5 & & \\
&               & & \properties_2 = -8 b_1 + 2 b_2 - 3 b_3 + 8 b_4 - 8 b_5 & & \\
&               & & \sum_{i=1}^5 b_i = 1 & & \\
&               & & \continuousvariables_1 \in [0, 16], \properties_1 \in [0, 20], \properties_2 \in [-10, 10], b \in \{0, 1\}^5.
\end{aligned}
\end{equation}
A binary variable $b_i$ is attached to the $i$th catalog item. Each property is a linear combination of all its possible catalog values with weights $b$. At the solution, $b_i$ is $1$ if item $i$ is selected, $0$ otherwise. The components of $b$ add up to $1$, that is a single catalog item is selected.

The optimal solution reported by SCIP\footnote{Since SCIP imposes that the objective function be linear, the model implements the epigraph reformulation.} is $(\continuousvariables_1, \properties_1, \properties_2, b_1, b_2, b_3, b_4, b_5) = (9.999999999, 3, 2, 0, 1, 0, 0, 0)$: the optimal catalog item is item $2$.
The SCIP optimization output is given in Table~\ref{tab:scip-scenario1}: the columns represent the CPU time in seconds, the number of branch-and-bound nodes, the number of simplex iterations, the number of rows in the current LP, the number of cuts added to the LP, the current dual (lower) bound, the current primal (upper) bound and the duality gap, respectively.
The solution is obtained after 4 iterations of the dual simplex algorithm, to which one cut has been added. This is a low computational cost, equivalent to that of the proposed approach.

\begin{table}[htbp!]
\caption{Scenario 1: SCIP optimization output.}
\label{tab:scip-scenario1}
\centering
\footnotesize
\begin{tabular}{cccccccc}
\hline
time (s)  & node   & simplex iterations & LP rows & cuts & dual bound & primal bound & gap \\
\hline
0.0  & 1      & 2         &  16 &   0  & 0.0 & 27.0 & $+\infty$ \\
0.0  & 1      & 3         &  16 &   0 & 16.0 & 27.0 & 68.75\% \\
0.1  & 1      & 3         &   1 &   0 & 16.0 & 27.0 & 68.75\% \\
0.1  & 1      & 4         &   2 &   1 & 27.0 & 27.0 & 0.00\% \\
0.1  & 1      & 4         &   2 &   1 & 27.0 & 27.0 & 0.00\% \\
\hline
\end{tabular}
\end{table}

\subsection{Scenario 2}
\label{sec:results-scenario2}

\subsubsection{Proposed approach}

\begin{figure}[htbp!]
\centering
\begin{subfigure}[b]{0.49\textwidth}
\includegraphics[width=\columnwidth]{example_scenario2}
\caption{Initial range of properties $\properties_1$ and $\properties_2$.}
\end{subfigure}
\hfill
\begin{subfigure}[b]{0.49\textwidth}
\includegraphics[width=\columnwidth]{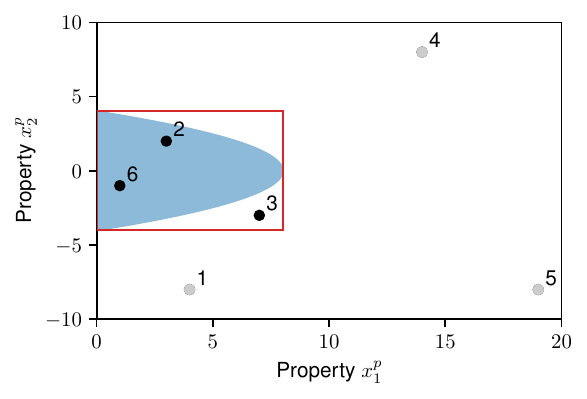}
\caption{HC4Revise with respect to $\constraints(\variables) = 0$.}
\end{subfigure} \\
\begin{subfigure}[b]{0.49\textwidth}
\includegraphics[width=\columnwidth]{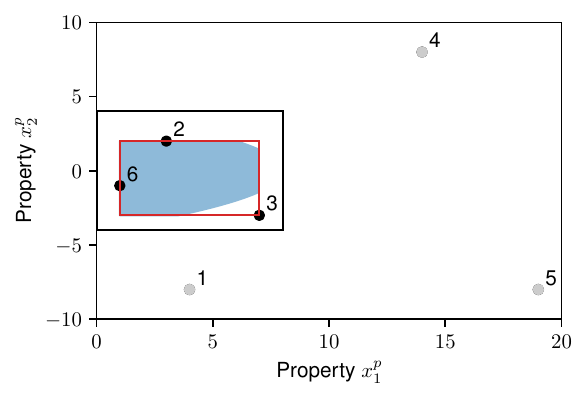}
\caption{\textsc{Clutch} with respect to $\catalogvariable(\properties) \in \catalog$.}
\end{subfigure}
\hfill
\begin{subfigure}[b]{0.49\textwidth}
\includegraphics[width=\columnwidth]{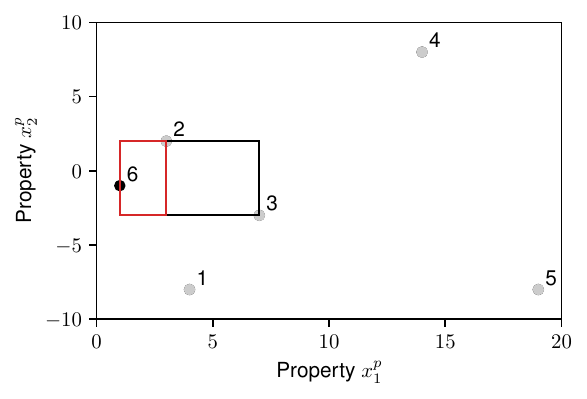}
\caption{HC4Revise with respect to $f \le \tilde{f} - \varepsilon$.}
\end{subfigure} \\
\begin{subfigure}[b]{0.49\textwidth}
\includegraphics[width=\columnwidth]{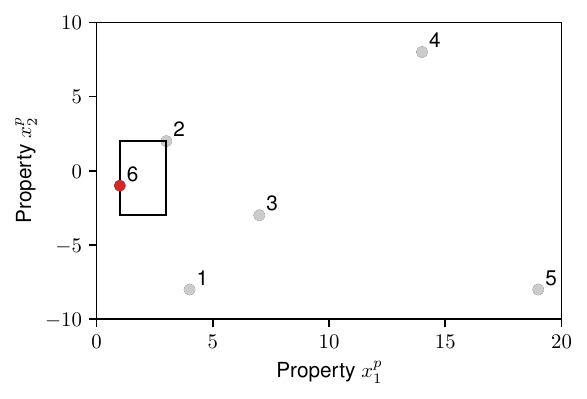}
\caption{\textsc{Clutch} with respect to $\catalogvariable(\properties) \in \catalog$.}
\end{subfigure}
\hfill
\begin{subfigure}[b]{0.49\textwidth}
\includegraphics[width=\columnwidth]{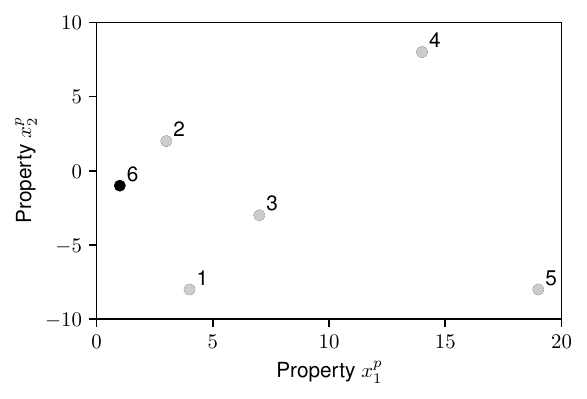}
\caption{HC4Revise with respect to $\constraints(\variables) = 0$.}
\end{subfigure}
\caption{Scenario 2: alternation of filtering and upper bounding phases.}
\label{fig:scenario2}
\end{figure}

Table~\ref{tab:scenario2-phases} describes the successive filtering and upper bounding phases that reduce the range of $x$ starting from $[0, 16] \times [0, 20] \times [-10, 10]$: the third column contains the current box for each phase and the fourth column presents the resulting ranges. The first column indicates the corresponding subfigures in Figure~\ref{fig:scenario2}. The problem was solved with an arbitrary tolerance $\varepsilon > 0$; the shortcut $3^-$ stands for $\sqrt[3]{27-\varepsilon}$, a quantity close to but strictly smaller than 3.

\begin{table}[htbp!]
\caption{Scenario 2: alternation of filtering and upper bounding phases.}
\label{tab:scenario2-phases}
\centering
\footnotesize
\begin{tabular}{clcc}
\hline
 & Phase                          & Current box                                & Result \\
\hline
(b) & HC4Revise($\constraints(\variables) = 0$)         & $[0, 16] \times [0, 20] \times [-10, 10]$  & $[0, 16] \times [0, 8] \times [-4, 4]$ \\
(c) & \textsc{Clutch}($\catalogvariable(\properties) \in \catalog$)   & $[0, 16] \times [0, 8] \times [-4, 4]$     & $[0, 16] \times [1, 7] \times [-3, 2]$ \\
& HC4Revise($\constraints(\variables) = 0$)             & $[0, 16] \times [1, 7] \times [-3, 2]$     & $[2, 16] \times [1, 7] \times [-3, 2]$ \\
(d) & HC4Revise($f(\variables) \le \tilde{f} - \varepsilon$)             & $[2, 16] \times [1, 7] \times [-3, 2]$     & $[2, 16] \times [1, 3^-] \times [-3, 2]$ \\
(e) & \textsc{Clutch}($\catalogvariable(\properties) \in \catalog$)   & $[2, 16] \times [1, 3^-] \times [-3, 2]$     & $[2, 16] \times [1, 1] \times [-1, -1]$ \\
(f) & HC4Revise($\constraints(\variables) = 0$)         & $[2, 16] \times [1, 1] \times [-1, -1]$      & $\mathbf{[3, 3] \times [1, 1] \times [-1, -1]}$ \\
\hline
\end{tabular}
\end{table}

Figure~\ref{fig:scenario2} provides a visual interpretation of the successive phases.
The first contraction step discards three of the catalog items. An alternation of HC4Revise and \textsc{Clutch} reduces the domain to the convex hull of catalog items 2, 3 and 6. An upper bounding phase (not shown in the figure) then selects the feasible item 2 (with properties $(\properties_1, \properties_2) = (3, 2)$) and the corresponding value $\continuousvariables_1 = 10$ and updates the best known upper bound: $\tilde{f} \gets f(10, 3, 2) = 27$. The subsequent contraction phase with respect to $f \le \tilde{f} - \varepsilon$ eliminates a large part of the domain (including item 2), which is then contracted to the optimal solution (item 6) after a \textsc{Clutch} phase. In this case, upper bounding and the use of the dynamic constraint $f(\variables) \le \tilde{f} - \varepsilon$ are required to discard one the two feasible items; however, no branching phase was required.

\subsubsection{Model with binary decision variables}

The corresponding model with binary decision variables is given by:
\begin{equation}
\begin{aligned}
& \underset{\variables, b}{\text{min}} & & (\properties_1)^3 & & \\
& \text{s.t.} 	& & \continuousvariables_1 - (\properties_2)^2 - 2 \properties_1 = 0 & & \\
&               & & \properties_1 = 4 b_1 + 3 b_2 + 7 b_3 + 14 b_4 + 19 b_5 + b_6 & & \\
&               & & \properties_2 = -8 b_1 + 2 b_2 - 3 b_3 + 8 b_4 - 8 b_5 - b_6 & & \\
&               & & \sum_{i=1}^6 b_i = 1 & & \\
&               & & \continuousvariables_1 \in [0, 16], \properties_1 \in [0, 20], \properties_2 \in [-10, 10], b \in \{0, 1\}^6
\end{aligned}
\end{equation}

The optimal solution reported by SCIP is $(\continuousvariables_1, \properties_1, \properties_2, b_1, b_2, b_3, b_4, b_5, b_6) = (3, 1, -1, 0, 0, 0, 0, 0, 1)$: the optimal catalog item is item $6$.
The SCIP optimization output is given in Table~\ref{tab:scip-scenario2}.
The solution is obtained after 24 iterations of the dual simplex algorithm, to which 3 cuts have been added. The computational cost is notably higher than that of SCIP for scenario 1 and that of the proposed approach for both scenarios.

\begin{table}[htbp!]
\caption{Scenario 2: SCIP optimization output.}
\label{tab:scip-scenario2}
\centering
\footnotesize
\begin{tabular}{ccccccccccc}
\hline
time (s) & node  & simplex iterations & LP rows & cuts & dual bound & primal bound &  gap \\  
\hline
0.0 &     1 &     2 &  19 &  0 & 0.0 & 27.0 &    $+\infty$ \\
0.0 &     1 &     9 &  29 &  0 & 0.0 & 27.0 &    $+\infty$ \\
0.0 &     1 &    13 &  31 &   2 & 0.0 & 27.0 &    $+\infty$ \\
0.0 &     1 &    15 &  32 &   3 & 1.0 & 27.0 & 2600.00\% \\
0.0 &     1 &    15 &  32 &   3 & 1.0 & 27.0 & 2600.00\% \\
0.0 &     1 &    24 &  33 &   3 & 1.0 & 27.0 & 2600.00\% \\
0.0 &     1 &    24 &  33 &   3 & 1.0 & 1.0 &   0.00\% \\
0.0 &     1 &    24 &  33 &   3 & 1.0 & 1.0 &   0.00\% \\
\hline
\end{tabular}
\end{table}

\section{Conclusion}

We described a rigorous interval constraint programming framework that interleaves branching and filtering phases, and treats the underlying catalog properties of the categorical variables as optimization variables. No additional variable is introduced and no modeling effort is required from the user.
We introduced \textsc{Clutch}, a novel catalog-based contractor invoked upon branching or filtering; it guarantees consistency between the categorical properties and the catalog items during the arborescent search. We demonstrated how \textsc{Clutch} can be easily implemented in off-the-shelf interval-based continuous solvers (such as Ibex, IBBA, GlobSol or Charibde) that solve the problem to global optimality, even in the presence of roundoff errors.
We introduced a toy problem with one continuous variable and two categorical properties and used the proposed approach to solve two scenarios to global optimality: it alternates between filtering (\textsc{Clutch} and classical contractors), upper and lower bounding, and branching to iteratively reduce the ranges of the continuous and property variables until convergence to the global minimum.
The toy problem can be solved efficiently with little computational power compared to an equivalent MINLP model with binary variables, while retaining the advantages of interval-based methods (reliability and maturity of existing codes).

\section*{Acknowledgments}
The author is grateful to the anonymous referees and the editor for their insightful comments that improved both the content and the presentation of this article.

\begin{appendices}
\section{Overestimation in models with binary decision variables: an example}
\label{sec:overestimation}

Consider a catalog variable whose properties $\properties_1$ and $\properties_2$ are available from a catalog $\catalog$ (Table~\ref{tab:overestimation-catalog}).

\begin{table}[htbp!]
\caption{Catalog $\catalog$.}
\label{tab:overestimation-catalog}
\centering
\begin{tabular}{c|cc}
\hline
Item   & Property $\properties_1$    & Property $\properties_2$ \\
\hline
1       & -1                & 0.1 \\
2 	    & 0.5               & 0.2 \\
3 	    & 0.6               & -0.2 \\
\hline
\end{tabular}
\end{table}
\end{appendices}

The images of the catalog items by the function $f(\properties) = (\properties_1)^2 - 5 \properties_2$ are given by: $\{f(-1, 0.1), f(0.5, 0.2), f(0.6, -0.2)\} = \{0.5, -0.75, 1.36\}$. We now compare enclosures of the range of $f$ for a model with continuous relaxations of $\properties$ (similarly to the proposed approach) and a model with binary decision variables.

\subsection{Continuous relaxations}

According to the catalog, the property $\properties_1$ (resp. $\properties_2$) lives in the interval $[-1, 0.6]$ (resp. $[-0.2, 0.2]$). An interval enclosure of the range of $f(\properties)$ is given by:
\begin{equation}
(\properties_1)^2 - 5 \properties_2 \in [-1, 0.6]^2 - 5 [-0.2, 0.2] = [0, 1] - [-1, 1] = [-1, 2]
\end{equation}

This is an overestimation of the true range. However, this is the best possible interval enclosure, since the variables $\properties_1$ and $\properties_2$ have a single occurrence in the expression of $f$ and $f$ is continuous~\cite{Moore1966Interval}.

\subsection{Binary decision variables}

We introduce the binary decision variables $(b_1, b_2, b_3) \in \{0, 1\}^3$ and express $\properties_1$ and $\properties_2$ as functions of $(b_1, b_2, b_3)$ and the catalog values:
\begin{eqnarray}
\properties_1 = -1 \times b_1 + 0.5 \times b_2 + 0.6 \times b_3 \\
\properties_2 = 0.1 \times b_1 + 0.2 \times b_2 - 0.2 \times b_3
\end{eqnarray}
$(b_1, b_2, b_3)$ are linked by the choice constraint $b_1 + b_2 + b_3 = 1$.
We now relax the integrality constraints and compute an interval enclosure of the range of $f(\properties)$:
\begin{equation}
\begin{aligned}
(\properties_1)^2 - 5 \properties_2   & = (-b_1 + 0.5 b_2 + 0.6 b_3)^2 - 5 (0.1 b_1 + 0.2 b_2 - 0.2 b_3) \\
                & \in (-[0, 1] + 0.5 [0, 1] + 0.6 [0, 1])^2 - 5(0.1 [0, 1] + 0.2 [0, 1] - 0.2 [0, 1]) \\
                & = [-1, 1.1]^2 - 5 [-0.2, 0.3] \\
                & = [0, 1.21] - [-1, 1.5] \\
                & = [-1.5, 2.21] \\
\end{aligned}
\end{equation}
This is also an overestimation of the true range that is slightly larger than that produced with continuous relaxations. This is due to the multiple occurrences of the binary variables in the expression of $f$, a phenomenon called \textit{dependency effect}~\cite{Moore1966Interval}.
A different enclosure is obtained by exploiting the choice constraint and factorizing; let $b_3 = 1 - b_1 - b_2$. Then:

\begin{equation}
\begin{aligned}
(\properties_1)^2 - 5 \properties_2   & = (-b_1 + 0.5 b_2 + 0.6 (1 - b_1 - b_2))^2 - 5 (0.1 b_1 + 0.2 b_2 - 0.2 (1 - b_1 - b_2)) \\
                & = (-0.7 b_1 - 0.1 b_2 + 0.6)^2 - 5 (0.3 b_1 + 0.4 b_2 - 0.2) \\
                & \in (-0.7 [0, 1] - 0.1 [0, 1] + 0.6)^2 - 5 (0.3 [0, 1] + 0.4 [0, 1] - 0.2) \\
                & = [-0.2, 0.6]^2 - 5 [-0.2, 0.5] \\
                & = [0, 0.36] - [-1, 2.5] \\
                & = [-2.5, 1.36] \\
\end{aligned}
\end{equation}

The upper bound is now optimal (it corresponds to picking item 3), while the lower bound was degraded.

\section*{Declarations}

\paragraph{Funding}
This work was supported by the research project UNSEEN funded by the German Federal Ministry for Economic Affairs and Climate Action (grant number FKZ 03EI1004C).

\paragraph{Conflicts of interest/Competing interests} None

\paragraph{Availability of data and material} None

\paragraph{Code availability}
A Julia prototype is available as open-source software under the MIT license at \url{https://github.com/cvanaret/CateGOrical.jl}~.

\bibliography{biblio}   


\begin{thebibliography}{42}
\ifx \bisbn   \undefined \def \bisbn  #1{ISBN #1}\fi
\ifx \binits  \undefined \def \binits#1{#1}\fi
\ifx \bauthor  \undefined \def \bauthor#1{#1}\fi
\ifx \batitle  \undefined \def \batitle#1{#1}\fi
\ifx \bjtitle  \undefined \def \bjtitle#1{#1}\fi
\ifx \bvolume  \undefined \def \bvolume#1{\textbf{#1}}\fi
\ifx \byear  \undefined \def \byear#1{#1}\fi
\ifx \bissue  \undefined \def \bissue#1{#1}\fi
\ifx \bfpage  \undefined \def \bfpage#1{#1}\fi
\ifx \blpage  \undefined \def \blpage #1{#1}\fi
\ifx \burl  \undefined \def \burl#1{\textsf{#1}}\fi
\ifx \doiurl  \undefined \def \doiurl#1{\url{https://doi.org/#1}}\fi
\ifx \betal  \undefined \def \betal{\textit{et al.}}\fi
\ifx \binstitute  \undefined \def \binstitute#1{#1}\fi
\ifx \binstitutionaled  \undefined \def \binstitutionaled#1{#1}\fi
\ifx \bctitle  \undefined \def \bctitle#1{#1}\fi
\ifx \beditor  \undefined \def \beditor#1{#1}\fi
\ifx \bpublisher  \undefined \def \bpublisher#1{#1}\fi
\ifx \bbtitle  \undefined \def \bbtitle#1{#1}\fi
\ifx \bedition  \undefined \def \bedition#1{#1}\fi
\ifx \bseriesno  \undefined \def \bseriesno#1{#1}\fi
\ifx \blocation  \undefined \def \blocation#1{#1}\fi
\ifx \bsertitle  \undefined \def \bsertitle#1{#1}\fi
\ifx \bsnm \undefined \def \bsnm#1{#1}\fi
\ifx \bsuffix \undefined \def \bsuffix#1{#1}\fi
\ifx \bparticle \undefined \def \bparticle#1{#1}\fi
\ifx \barticle \undefined \def \barticle#1{#1}\fi
\bibcommenthead
\ifx \bconfdate \undefined \def \bconfdate #1{#1}\fi
\ifx \botherref \undefined \def \botherref #1{#1}\fi
\ifx \url \undefined \def \url#1{\textsf{#1}}\fi
\ifx \bchapter \undefined \def \bchapter#1{#1}\fi
\ifx \bbook \undefined \def \bbook#1{#1}\fi
\ifx \bcomment \undefined \def \bcomment#1{#1}\fi
\ifx \oauthor \undefined \def \oauthor#1{#1}\fi
\ifx \citeauthoryear \undefined \def \citeauthoryear#1{#1}\fi
\ifx \endbibitem  \undefined \def \endbibitem {}\fi
\ifx \bconflocation  \undefined \def \bconflocation#1{#1}\fi
\ifx \arxivurl  \undefined \def \arxivurl#1{\textsf{#1}}\fi
\csname PreBibitemsHook\endcsname

\bibitem{gao2018categorical}
\begin{barticle}
\bauthor{\bsnm{Gao}, \binits{H.}},
\bauthor{\bsnm{Breitkopf}, \binits{P.}},
\bauthor{\bsnm{Coelho}, \binits{R.F.}},
\bauthor{\bsnm{Xiao}, \binits{M.}}:
\batitle{Categorical structural optimization using discrete manifold learning
  approach and custom-built evolutionary operators}.
\bjtitle{Structural and Multidisciplinary Optimization}
\bvolume{58}(\bissue{1}),
\bfpage{215}--\blpage{228}
(\byear{2018})
\end{barticle}
\endbibitem

\bibitem{nedvelkova2019splitting}
\begin{barticle}
\bauthor{\bsnm{Ned{\v{e}}lkov{\'a}}, \binits{Z.}},
\bauthor{\bsnm{Cromvik}, \binits{C.}},
\bauthor{\bsnm{Lindroth}, \binits{P.}},
\bauthor{\bsnm{Patriksson}, \binits{M.}},
\bauthor{\bsnm{Str{\"o}mberg}, \binits{A.-B.}}:
\batitle{A splitting algorithm for simulation-based optimization problems with
  categorical variables}.
\bjtitle{Engineering Optimization}
\bvolume{51}(\bissue{5}),
\bfpage{815}--\blpage{831}
(\byear{2019})
\end{barticle}
\endbibitem

\bibitem{barjhoux2017mixed}
\begin{bchapter}
\bauthor{\bsnm{Barjhoux}, \binits{P.-J.}},
\bauthor{\bsnm{Diouane}, \binits{Y.}},
\bauthor{\bsnm{Grihon}, \binits{S.}},
\bauthor{\bsnm{Bettebghor}, \binits{D.}},
\bauthor{\bsnm{Morlier}, \binits{J.}}:
\bctitle{Mixed variable structural optimization: toward an efficient hybrid
  algorithm}.
In: \bbtitle{World Congress of Structural and Multidisciplinary Optimisation},
pp. \bfpage{1880}--\blpage{1896}
(\byear{2017}).
\bcomment{Springer}
\end{bchapter}
\endbibitem

\bibitem{lindroth2011pure}
\begin{bbook}
\bauthor{\bsnm{Lindroth}, \binits{P.}},
\bauthor{\bsnm{Patriksson}, \binits{M.}}:
\bbtitle{Pure Categorical Optimization: a Global Descent Approach}.
\bpublisher{Department of Mathematical Sciences, Division of Mathematics,
  Chalmers}, \blocation{}
(\byear{2011})
\end{bbook}
\endbibitem

\bibitem{abhishek2010modeling}
\begin{barticle}
\bauthor{\bsnm{Abhishek}, \binits{K.}},
\bauthor{\bsnm{Leyffer}, \binits{S.}},
\bauthor{\bsnm{Linderoth}, \binits{J.T.}}:
\batitle{Modeling without categorical variables: a mixed-integer nonlinear
  program for the optimization of thermal insulation systems}.
\bjtitle{Optimization and Engineering}
\bvolume{11}(\bissue{2}),
\bfpage{185}--\blpage{212}
(\byear{2010})
\end{barticle}
\endbibitem

\bibitem{agurok2019multi}
\begin{botherref}
\oauthor{\bsnm{Agurok}, \binits{I.}}:
Multi-extremum optimization in lens design: navigation through merit function
  valleys maze.
arXiv preprint arXiv:1907.08676
(2019)
\end{botherref}
\endbibitem

\bibitem{lin2017advanced}
\begin{bbook}
\bauthor{\bsnm{Lin}, \binits{P.D.}}:
\bbtitle{Advanced Geometrical Optics}.
\bpublisher{Springer}, \blocation{}
(\byear{2017})
\end{bbook}
\endbibitem

\bibitem{lecoutre2012path}
\begin{bchapter}
\bauthor{\bsnm{Lecoutre}, \binits{C.}},
\bauthor{\bsnm{Likitvivatanavong}, \binits{C.}},
\bauthor{\bsnm{Yap}, \binits{R.H.}}:
\bctitle{A path-optimal {GAC} algorithm for table constraints}.
In: \bbtitle{ECAI 2012},
pp. \bfpage{510}--\blpage{515}.
\bpublisher{IOS Press}, \blocation{}
(\byear{2012})
\end{bchapter}
\endbibitem

\bibitem{brown1993solving}
\begin{barticle}
\bauthor{\bsnm{Brown}, \binits{D.R.}},
\bauthor{\bsnm{Hwang}, \binits{K.-Y.}}:
\batitle{Solving fixed configuration problems with genetic search}.
\bjtitle{Research in Engineering Design}
\bvolume{5}(\bissue{2}),
\bfpage{80}--\blpage{87}
(\byear{1993})
\end{barticle}
\endbibitem

\bibitem{carlson1996genetic}
\begin{barticle}
\bauthor{\bsnm{Carlson}, \binits{S.E.}}:
\batitle{Genetic algorithm attributes for component selection}.
\bjtitle{Research in Engineering Design}
\bvolume{8}(\bissue{1}),
\bfpage{33}--\blpage{51}
(\byear{1996})
\end{barticle}
\endbibitem

\bibitem{carlson1998using}
\begin{barticle}
\bauthor{\bsnm{Carlson-Skalak}, \binits{S.}},
\bauthor{\bsnm{White}, \binits{M.D.}},
\bauthor{\bsnm{Teng}, \binits{Y.}}:
\batitle{Using an evolutionary algorithm for catalog design}.
\bjtitle{Research in Engineering Design}
\bvolume{10}(\bissue{2}),
\bfpage{63}--\blpage{83}
(\byear{1998})
\end{barticle}
\endbibitem

\bibitem{neumaier2005comparison}
\begin{barticle}
\bauthor{\bsnm{Neumaier}, \binits{A.}},
\bauthor{\bsnm{Shcherbina}, \binits{O.}},
\bauthor{\bsnm{Huyer}, \binits{W.}},
\bauthor{\bsnm{Vink{\'o}}, \binits{T.}}:
\batitle{A comparison of complete global optimization solvers}.
\bjtitle{Mathematical programming}
\bvolume{103}(\bissue{2}),
\bfpage{335}--\blpage{356}
(\byear{2005})
\end{barticle}
\endbibitem

\bibitem{sahinidis1996baron}
\begin{barticle}
\bauthor{\bsnm{Sahinidis}, \binits{N.V.}}:
\batitle{Baron: A general purpose global optimization software package}.
\bjtitle{Journal of global optimization}
\bvolume{8}(\bissue{2}),
\bfpage{201}--\blpage{205}
(\byear{1996})
\end{barticle}
\endbibitem

\bibitem{Moore1966Interval}
\begin{bbook}
\bauthor{\bsnm{Moore}, \binits{R.E.}}:
\bbtitle{Interval Analysis}
vol. \bseriesno{4}.
\bpublisher{Prentice-Hall Englewood Cliffs}, \blocation{}
(\byear{1966})
\end{bbook}
\endbibitem

\bibitem{Trombettoni2011Inner}
\begin{bchapter}
\bauthor{\bsnm{Trombettoni}, \binits{G.}},
\bauthor{\bsnm{Araya}, \binits{I.}},
\bauthor{\bsnm{Neveu}, \binits{B.}},
\bauthor{\bsnm{Chabert}, \binits{G.}}:
\bctitle{Inner regions and interval linearizations for global optimization.}
In: \bbtitle{AAAI}
(\byear{2011})
\end{bchapter}
\endbibitem

\bibitem{Ninin2010Reliable}
\begin{bbook}
\bauthor{\bsnm{Ninin}, \binits{J.}},
\bauthor{\bsnm{Hansen}, \binits{P.}},
\bauthor{\bsnm{Messine}, \binits{F.}}:
\bbtitle{A Reliable Affine Relaxation Method for Global Optimization}.
\bpublisher{Groupe d'{\'e}tudes et de recherche en analyse des d{\'e}cisions},
  \blocation{}
(\byear{2010})
\end{bbook}
\endbibitem

\bibitem{Kearfott1996Rigorous}
\begin{bbook}
\bauthor{\bsnm{Kearfott}, \binits{R.B.}}:
\bbtitle{Rigorous Global Search: Continuous Problems}.
\bpublisher{Springer}, \blocation{}
(\byear{1996})
\end{bbook}
\endbibitem

\bibitem{Vanaret2015Hybridization}
\begin{botherref}
\oauthor{\bsnm{Vanaret}, \binits{C.}}:
Hybridization of interval methods and evolutionary algorithms for solving
  difficult optimization problems.
PhD thesis,
INP Toulouse
(2015)
\end{botherref}
\endbibitem

\bibitem{neumaier2004complete}
\begin{barticle}
\bauthor{\bsnm{Neumaier}, \binits{A.}}:
\batitle{Complete search in continuous global optimization and constraint
  satisfaction}.
\bjtitle{Acta numerica}
\bvolume{13},
\bfpage{271}--\blpage{369}
(\byear{2004})
\end{barticle}
\endbibitem

\bibitem{Puget1994ilogsolver}
\begin{bchapter}
\bauthor{\bsnm{Puget}, \binits{J.-F.}}:
\bctitle{A {C++} implementation of {CLP}}.
In: \bbtitle{Procs. of the Singapore Conference on Intelligent Systems}
(\byear{1994})
\end{bchapter}
\endbibitem

\bibitem{van1997numerica}
\begin{bbook}
\bauthor{\bsnm{Van~Hentenryck}, \binits{P.}},
\bauthor{\bsnm{Michel}, \binits{L.}},
\bauthor{\bsnm{Deville}, \binits{Y.}}:
\bbtitle{Numerica: a Modeling Language for Global Optimization}.
\bpublisher{MIT press}, \blocation{}
(\byear{1997})
\end{bbook}
\endbibitem

\bibitem{jussien2008choco}
\begin{bchapter}
\bauthor{\bsnm{Jussien}, \binits{N.}},
\bauthor{\bsnm{Rochart}, \binits{G.}},
\bauthor{\bsnm{Lorca}, \binits{X.}}:
\bctitle{Choco: an open source java constraint programming library}.
In: \bbtitle{CPAIOR'08 Workshop on Open-Source Software for Integer and
  Contraint Programming (OSSICP'08)},
pp. \bfpage{1}--\blpage{10}
(\byear{2008})
\end{bchapter}
\endbibitem

\bibitem{chabert2009contractor}
\begin{barticle}
\bauthor{\bsnm{Chabert}, \binits{G.}},
\bauthor{\bsnm{Jaulin}, \binits{L.}}:
\batitle{Contractor programming}.
\bjtitle{Artificial Intelligence}
\bvolume{173}(\bissue{11}),
\bfpage{1079}--\blpage{1100}
(\byear{2009})
\end{barticle}
\endbibitem

\bibitem{Benhamou1999Revising}
\begin{bchapter}
\bauthor{\bsnm{Benhamou}, \binits{F.}},
\bauthor{\bsnm{Goualard}, \binits{F.}},
\bauthor{\bsnm{Granvilliers}, \binits{L.}},
\bauthor{\bsnm{Puget}, \binits{J.-F.}}:
\bctitle{Revising hull and box consistency}.
In: \bbtitle{International Conference on Logic Programming},
pp. \bfpage{230}--\blpage{244}.
\bpublisher{MIT press}, \blocation{}
(\byear{1999})
\end{bchapter}
\endbibitem

\bibitem{zhang2007new}
\begin{bchapter}
\bauthor{\bsnm{Zhang}, \binits{X.}},
\bauthor{\bsnm{Liu}, \binits{S.}}:
\bctitle{A new interval-genetic algorithm}.
In: \bbtitle{Third International Conference on Natural Computation (ICNC
  2007)},
vol. \bseriesno{4},
pp. \bfpage{193}--\blpage{197}
(\byear{2007}).
\bcomment{IEEE}
\end{bchapter}
\endbibitem

\bibitem{gallardo2007hybridization}
\begin{barticle}
\bauthor{\bsnm{Gallardo}, \binits{J.E.}},
\bauthor{\bsnm{Cotta}, \binits{C.}},
\bauthor{\bsnm{Fern{\'a}ndez}, \binits{A.J.}}:
\batitle{On the hybridization of memetic algorithms with branch-and-bound
  techniques}.
\bjtitle{IEEE Transactions on Systems, Man, and Cybernetics, Part B
  (Cybernetics)}
\bvolume{37}(\bissue{1}),
\bfpage{77}--\blpage{83}
(\byear{2007})
\end{barticle}
\endbibitem

\bibitem{blum2011hybrid}
\begin{barticle}
\bauthor{\bsnm{Blum}, \binits{C.}},
\bauthor{\bsnm{Puchinger}, \binits{J.}},
\bauthor{\bsnm{Raidl}, \binits{G.R.}},
\bauthor{\bsnm{Roli}, \binits{A.}}:
\batitle{Hybrid metaheuristics in combinatorial optimization: A survey}.
\bjtitle{Applied soft computing}
\bvolume{11}(\bissue{6}),
\bfpage{4135}--\blpage{4151}
(\byear{2011})
\end{barticle}
\endbibitem

\bibitem{cotta1995hybridizing}
\begin{bchapter}
\bauthor{\bsnm{Cotta}, \binits{C.}},
\bauthor{\bsnm{Aldana}, \binits{J.}},
\bauthor{\bsnm{Nebro}, \binits{A.J.}},
\bauthor{\bsnm{Troya}, \binits{J.M.}}:
\bctitle{Hybridizing genetic algorithms with branch and bound techniques for
  the resolution of the {TSP}}.
In: \bbtitle{Artificial Neural Nets and Genetic Algorithms},
pp. \bfpage{277}--\blpage{280}
(\byear{1995}).
\bcomment{Springer}
\end{bchapter}
\endbibitem

\bibitem{alliot2012finding}
\begin{bchapter}
\bauthor{\bsnm{Alliot}, \binits{J.-M.}},
\bauthor{\bsnm{Durand}, \binits{N.}},
\bauthor{\bsnm{Gianazza}, \binits{D.}},
\bauthor{\bsnm{Gotteland}, \binits{J.-B.}}:
\bctitle{Finding and proving the optimum: Cooperative stochastic and
  deterministic search.}
In: \bbtitle{ECAI},
pp. \bfpage{55}--\blpage{60}
(\byear{2012})
\end{bchapter}
\endbibitem

\bibitem{ichida1979interval}
\begin{barticle}
\bauthor{\bsnm{Ichida}, \binits{K.}},
\bauthor{\bsnm{Fujii}, \binits{Y.}}:
\batitle{An interval arithmetic method for global optimization}.
\bjtitle{Computing}
\bvolume{23}(\bissue{1}),
\bfpage{85}--\blpage{97}
(\byear{1979})
\end{barticle}
\endbibitem

\bibitem{kearfott1990algorithm}
\begin{barticle}
\bauthor{\bsnm{Kearfott}, \binits{R.B.}},
\bauthor{\bsnm{Novoa~III}, \binits{M.}}:
\batitle{Algorithm 681: {INTBIS}, a portable interval {N}ewton/bisection
  package}.
\bjtitle{ACM Transactions on Mathematical Software (TOMS)}
\bvolume{16}(\bissue{2}),
\bfpage{152}--\blpage{157}
(\byear{1990})
\end{barticle}
\endbibitem

\bibitem{sherali2013reformulation}
\begin{bbook}
\bauthor{\bsnm{Sherali}, \binits{H.D.}},
\bauthor{\bsnm{Adams}, \binits{W.P.}}:
\bbtitle{A Reformulation-linearization Technique for Solving Discrete and
  Continuous Nonconvex Problems}
vol. \bseriesno{31}.
\bpublisher{Springer}, \blocation{}
(\byear{2013})
\end{bbook}
\endbibitem

\bibitem{Araya2012Contractor}
\begin{bchapter}
\bauthor{\bsnm{Araya}, \binits{I.}},
\bauthor{\bsnm{Trombettoni}, \binits{G.}},
\bauthor{\bsnm{Neveu}, \binits{B.}}:
\bctitle{A contractor based on convex interval {T}aylor}.
In: \bbtitle{International Conference on Integration of Artificial Intelligence
  (AI) and Operations Research (OR) Techniques in Constraint Programming},
pp. \bfpage{1}--\blpage{16}
(\byear{2012}).
\bcomment{Springer}
\end{bchapter}
\endbibitem

\bibitem{sunaga1958theory}
\begin{barticle}
\bauthor{\bsnm{Sunaga}, \binits{T.}}:
\batitle{Theory of interval algebra and its application to numerical analysis}.
\bjtitle{RAAG memoirs}
\bvolume{2}(\bissue{29-46}),
\bfpage{209}
(\byear{1958})
\end{barticle}
\endbibitem

\bibitem{Neumaier2004Safe}
\begin{barticle}
\bauthor{\bsnm{Neumaier}, \binits{A.}},
\bauthor{\bsnm{Shcherbina}, \binits{O.}}:
\batitle{Safe bounds in linear and mixed-integer linear programming}.
\bjtitle{Mathematical Programming}
\bvolume{99}(\bissue{2}),
\bfpage{283}--\blpage{296}
(\byear{2004})
\end{barticle}
\endbibitem

\bibitem{messine1997methodes}
\begin{botherref}
\oauthor{\bsnm{Messine}, \binits{F.}}:
M{\'e}thodes d'optimisation globale bas{\'e}es sur l'analyse d'intervalle pour
  la r{\'e}solution de probl{\`e}mes avec contraintes.
PhD thesis,
Toulouse, INPT
(1997)
\end{botherref}
\endbibitem

\bibitem{messine2004deterministic}
\begin{barticle}
\bauthor{\bsnm{Messine}, \binits{F.}}:
\batitle{Deterministic global optimization using interval constraint
  propagation techniques}.
\bjtitle{RAIRO-Operations Research-Recherche Op{\'e}rationnelle}
\bvolume{38}(\bissue{4}),
\bfpage{277}--\blpage{293}
(\byear{2004})
\end{barticle}
\endbibitem

\bibitem{belotti2010feasibility}
\begin{bchapter}
\bauthor{\bsnm{Belotti}, \binits{P.}},
\bauthor{\bsnm{Cafieri}, \binits{S.}},
\bauthor{\bsnm{Lee}, \binits{J.}},
\bauthor{\bsnm{Liberti}, \binits{L.}}:
\bctitle{Feasibility-based bounds tightening via fixed points}.
In: \bbtitle{International Conference on Combinatorial Optimization and
  Applications},
pp. \bfpage{65}--\blpage{76}
(\byear{2010}).
\bcomment{Springer}
\end{bchapter}
\endbibitem

\bibitem{cleary1987logical}
\begin{barticle}
\bauthor{\bsnm{Cleary}, \binits{J.G.}}:
\batitle{Logical arithmetic}.
\bjtitle{Future computing systems}
\bvolume{2}(\bissue{2}),
\bfpage{125}--\blpage{149}
(\byear{1987})
\end{barticle}
\endbibitem

\bibitem{mackworth1977consistency}
\begin{barticle}
\bauthor{\bsnm{Mackworth}, \binits{A.K.}}:
\batitle{Consistency in networks of relations}.
\bjtitle{Artificial intelligence}
\bvolume{8}(\bissue{1}),
\bfpage{99}--\blpage{118}
(\byear{1977})
\end{barticle}
\endbibitem

\bibitem{zamora1999branch}
\begin{barticle}
\bauthor{\bsnm{Zamora}, \binits{J.M.}},
\bauthor{\bsnm{Grossmann}, \binits{I.E.}}:
\batitle{A branch and contract algorithm for problems with concave univariate,
  bilinear and linear fractional terms}.
\bjtitle{Journal of Global Optimization}
\bvolume{14},
\bfpage{217}--\blpage{249}
(\byear{1999})
\end{barticle}
\endbibitem

\bibitem{SCIP}
\begin{botherref}
\oauthor{\bsnm{Bestuzheva}, \binits{K.}}, et al.:
{The SCIP Optimization Suite 8.0}.
ZIB-Report 21-41,
Zuse Institute Berlin
(December 2021)
\end{botherref}
\endbibitem

\end{thebibliography}

\end{document}